\documentclass[11pt]{amsart}
\usepackage{amsmath,paralist}
\usepackage{amsthm}
\usepackage{graphicx}
\usepackage{tikz,color}
\usepackage[margin=1in]{geometry}
\usetikzlibrary{arrows,shapes,snakes,backgrounds} 
\usepackage{hyperref}
 
\newtheorem{theorem}{Theorem}

\newtheorem{lemma}[theorem]{Lemma}
\newtheorem{corollary}[theorem]{Corollary}

\theoremstyle{remark}

\theoremstyle{definition}

\def\ca{{\rm Cox^A(n)}}
\def\cb{{\rm Cox^{BC}(n)}}
\def\cc{{\rm Cox^{BC}(n)}}

\def\x{{\bf x}}
\def\l{{\lambda}}
 
\def\c{{\bf c}}
\def\o{{\bf o }}

\def\a{{\alpha}}
 \def\f_H{{\bf w}}
 
 \def\a{{\bf a}}
\def\R{\mathbb{R}}

\def\aa{{\bf a}}
 \def\A{\mathcal{A}}
 \def\H{\mathcal{H}}
 \def\S{\mathfrak{S}}
 \def\sa{\mathcal{S}^A_n}
  
   \def\scc{\mathcal{S}^C_n}
\def\san{{{(\mathcal{S}^A_n)}^{\sqcup (n+1)}}}

\begin{document}

\tikzstyle{w}=[label=right:$\textcolor{red}{\cdots}$] 
\tikzstyle{b}=[label=right:$\cdot\,\textcolor{red}{\cdot}\,\cdot$] 
\tikzstyle{bb}=[circle,draw=black!90,fill=black!100,thick,inner sep=2pt,minimum width=3pt] 
\tikzstyle{b2}=[label=right:$\cdots$] 
\tikzstyle{w2}=[]
\tikzstyle{vw}=[label=above:$\textcolor{red}{\vdots}$] 
\tikzstyle{vb}=[label=above:$\vdots$] 

\tikzstyle{level 1}=[level distance=3.5cm, sibling distance=3.5cm]
\tikzstyle{level 2}=[level distance=3.5cm, sibling distance=2cm]

\tikzstyle{bag} = [text width=4em, text centered]
\tikzstyle{end} = [circle, minimum width=3pt,fill, inner sep=0pt]

\title[Posets, parking functions and the regions of the Shi arrangement]{Posets, parking functions and the regions of the Shi arrangement revisited}
\author{Karola M\'esz\'aros}
\address{
Department of Mathematics, Massachusetts Institute of Technology, Cambridge, MA 02139
}
\begin{abstract}
  The number of regions of the type $A_{n-1}$ Shi arrangement in $\R^n$ is counted by the 
  intrinsically   beautiful  formula $(n+1)^{n-1}$. First proved by Shi, this result motivated Pak and  Stanley as well as Athanasiadis and Linusson to provide bijective proofs. We give a description of the Athanasiadis-Linusson bijection and generalize it to a bijection between the regions of the  type $C_n$  Shi arrangement  in $\R^n$ and sequences $a_1a_2 \ldots a_n$, where $a_i \in \{-n, -n+1, \ldots, -1, 0, 1, \ldots, n-1, n\}$, $ i \in [n]$.  
Our bijections naturally restrict to bijections between regions  of the arrangements with a certain number of ceilings (or floors) and sequences with a given number of distinct elements. A special family of posets, whose antichains encode the regions of the arrangements, play a central role in our approach.
  
  \end{abstract}

\maketitle
  
  \section{Introduction}
  \label{sec:intro}
    
    A \textbf{hyperplane arrangement}  $\A$ is a finite set of affine hyperplanes in $\R^n$. The \textbf{regions} of $\A$ are the connected components of the space $\R^n \backslash \cup_{H \in \A} H$. 
       In this paper we  study \textbf{Shi arrangements}  of type $A_{n-1}$ and $C_n$, which are affine  hyperplane arrangements whose hyperplanes are parallel to reflecting hyperplanes of Coxeter groups. Denote by ${\rm Cox}^A(n)$ the \textbf{Coxeter arrangement of type $A_{n-1}$:}
    
   $${\rm Cox}^A(n)=\{x_i-x_j=0 \mid 1\leq i<j\leq n\}.$$
    
    Note that the regions of  ${\rm Cox}^A(n)$ are naturally indexed by type $A_{n-1}$ permutations $w \in \mathfrak{S}_n$. Namely, if $C^A\subset \R^n$ is the \textbf{dominant cone} of  ${\rm Cox}^A(n)$  defined by $x_1>x_2>\cdots>x_n$, then 
    
    $$ wC^A=\{ \x \in \R^n | x_{w(1)}>x_{w(2)}>\cdots >x_{w(n)} \}.$$ 
    
    Thus, the number of regions of $\ca$ is $|\mathfrak{S}_n|=n!.$ The \textbf{type $A_{n-1}$ Shi arrangement $\sa$} was first defined by Shi \cite{shi}:
    
    $$\sa=\ca \cup \{x_i-x_j=1 \mid 1\leq i<j\leq n\}.$$
    
    Shi \cite{shi}  proved the beautiful result that the number of regions of  $\mathcal{S}_n^A$ is $(n+1)^{n-1}$.  This statement is clearly  deserving of a combinatorial proof; two different bijections proving this result were provided by 
   Stanley \cite{sta0, sta} and  Athanasiadis and Linusson \cite{ath-lin}.  We give a  description of the Athanasiadis-Linusson bijection and generalize it to  type $C_n$, thereby answering a question of Athanasiadis and Linusson
   \cite[Section 4, Question 3]{ath-lin}. We also study statistics naturally arising from the bijections. In their forthcoming work on parking spaces \cite{arm-rei-rho}, Armstrong, Reiner and Rhoades provide the ultimate generalization of the Athanasiadis-Linusson bijection by constructing a uniform bijection for all crystallographic root systems.
    
    \medskip
     
 We now review the definitions necessary to state our results.
 
  A sequence $\a=(a_1, a_2, \ldots, a_n) \in [n]^n$ is a \textbf{parking function} if and only if 
 the increasing rearrangement $b_1\leq b_2\leq \cdots \leq b_n$ of $a_1, a_2, \ldots, a_n$ satisfies $b_i \leq i$. Denote by $PF(n)$ the set of all parking functions of length $n$. 
  Let $$\A(n)=\{\a=(a_1, a_2, \ldots, a_n) | a_i \in [n+1], i \in [n]\}.$$ Denote by $d(\a)$ the number of distinct numbers contained in the sequence $\a$.
  
  A hyperplane $H$ is a \textbf{wall} of a region $R$ if it is the affine span of a codimension-1 face of $R$. A wall $H$ is called a \textbf{floor} if $H$ does not contain the origin and $R$ and the origin lie in opposite half-spaces defined by $H$. Denote by $f(R)$ the number of floors of $R$. A wall $H$ is called a  \textbf{ceiling} if $H$ does not contain the origin and $R$ and the origin lie in the same half-spaces defined by $H$. Denote by $c(R)$ the number of floors of $R$. Denote by $R(\mathcal{H})$ the set of regions of the hyperplane arrangement $\mathcal{H}$.
   
   Our first description of the type $A_{n-1}$ bijection proving that $|R(\sa)|=(n+1)^{n-1}$ yields a natural correspondence between the multiset $\mathcal{M}$ in which each element of $R(\sa)$ appears $n+1$ times and the sequences $\A(n)$, while the second description, analogous to that of Athanasiadis and Linusson \cite{ath-lin}  is a direct correspondence between $R(\sa)$ and parking functions. The properties of these bijections yield Theorem \ref{gf}.
 
\begin{theorem} \cite{ath-lin} There is a bijection between the regions of $\sa$ and parking functions of length $n$.

\end{theorem}

 \begin{theorem} \label{gf} 
$$ \sum_{R \in R(\sa)} q^{c(R)}= \sum_{R \in R(\sa)} q^{f(R)}=\frac{1}{n+1}\sum_{\textrm{a} \in \A(n)} q^{n-d(\a)}=\sum_{\a \in PF(n)} q^{n-d(\a)}.$$
  \end{theorem}
 
 We use the techniques developed  for the type $A_{n-1}$ case to construct bijective proofs for type $C_n$. 
 
 The  \textbf{type $B_{n}$ and $C_n$  Coxeter arrangement $\cb$} in $\R^n$ is defined as follows. 
\begin{align*} \cb&=\{x_i-x_j=0, x_i+x_j=0, x_k=0 \mid 1\leq i<j\leq n, k \in [n]\}, \\     &=\{x_i-x_j=0, x_i+x_j=0, 2x_k=0 \mid 1\leq i<j\leq n, k \in [n]\}.\end{align*}

    Just as in the type $A_{n-1}$ case, the regions of the arrangements $\cb$ naturally correspond to type $B_n$ permutations $w \in \S_n^B$.   Recall that $\S_n^B$ is the group  of all bijections $w$ of the set $[\pm n]=\{-n, -n+1, \ldots, -1, 1, \ldots, n-1, n\}$ onto itself such that $$w(-i)=-w(i),$$ for all $i \in [\pm n]$ and composition as group operation. The notation $w=[a_1, \ldots, a_n]$ means $w(i)=a_i$, for $i \in [n]$, and is called the \textbf{window} of $w$.

     Let $C^{BC}\subset \R^n$ be the \textbf{dominant cone} of  $\cc$  defined by $$C^{BC}=\{\x \in \R^n \mid -x_n>-x_{n-1}>\cdots>-x_2>-x_1>x_1>x_2>\cdots>x_{n-1}>x_n\}.$$ 
    
    Let 
    $$ wC^{BC}=\{ \x \in \R^n | x_{w(-n)}>x_{w(-n+1)}>\cdots >x_{w(-1)}> x_{w(1)}>x_{w(2)}>\cdots >x_{w(n)} \},$$
    where $\{x_1, \ldots, x_n\}$ are the standard coordinate functions on $\R^n$ and $x_{-i}=-x_i$ for $i <0$. It follows that the number of regions of $\cb$ is $|\mathfrak{S}^B_n|=2^n n!.$ 

    \medskip
    
The      \textbf{type $C_n$  Shi arrangement  $\scc$}  is as expected: 
     
     $$\scc=\cc \cup \{x_i-x_j=1, x_i+x_j=1, 2x_k=1 \mid 1\leq i<j\leq n, k \in [n]\}.$$

We construct bijections between the  regions of the  type $C_{n}$  Shi arrangement $\scc$ in $\R^n$ and sequences in the set $$\A^{C}(n)=\{ (a_1, a_2, \ldots, a_n)| a_i \in \{-n, -n+1, \ldots, -1, 0, 1, \ldots, n-1, n\},  i \in [n]\}.$$  
Athanasiadis and Linusson \cite[Section 4, Question 3]{ath-lin} were the first to ask for the construction of such  bijection in their paper dealing with the type $A_{n-1}$ case. The properties of our bijections yield Theorem \ref{bgf}.

\begin{theorem} \label{thm:bijc} There is a bijection between the regions of $\scc$ and  sequences in the set $\A^{C}(n)$.
\end{theorem}

 \begin{theorem} \label{bgf} 
$$ \sum_{R \in R(\scc)} q^{c(R)}= \sum_{R \in R(\scc)} q^{f(R)}=\sum_{\a \in \A^{BC}(n)} q^{n-d^{C}(\a)},$$

 \noindent where  $d^{C}(\a)$ is the number of distinct absolute values of the nonzero numbers appearing in ~$\a$.
 \end{theorem}

The outline of the paper is as follows. In Section \ref{sec:poset} we explain the connection between the regions of Shi arrangements and posets of nonnesting partitions. In  Section \ref{secaa} we build on this connection to give a description of the Athanasiadis-Linusson bijection between the regions of the type $A_{n-1}$ Shi arrangement and parking functions, as well as prove  Theorem \ref{gf}. In Section \ref{secc} we generalize the contents of Section \ref{secaa} to the type  $C_n$ case,  proving Theorems \ref{thm:bijc} and \ref{bgf}.  Section \ref{sec:aa} is the story of  Section \ref{secaa} without  arrangements, only in terms of posets and sequences of type $A_{n-1}$. Section  \ref{sec:bb} similarly reiterates the basic thoughts in type $C_n$ on the level of posets and sequences. 

        \section{Posets and the regions of Shi arrangements}
    \label{sec:poset}
    
    Our bijections  are based on a correspondence developed by Stanley in \cite[Section 5]{sta0}  between the antichains of a certain family of posets  $Q_w$, $w \in \S_n$, and regions of $\sa$. In this section we explain this correspondence and its type $C_n$ extension.  
    For a related  bijection  between the positive chambers of the Shi arrangement and order ideals of the root poset of  corresponding type  see \cite[Theorem 5.1.13]{arm} and \cite{cel-pap}.
    For basic definitions about posets  see \cite[Chapter 3]{ec1}.
    
   \subsection{Poset and the regions of $\sa$.}   \label{sub:aa}  Each region of $\sa$ lies in one of the cones  of ${\rm Cox}^A(n)$. We restrict our attention to the regions of $\sa$ in an arbitrary cone $wC^A$, $w \in \S_n$. Each such region is uniquely determined by the set of its ceilings (or the set of its floors). The set of hyperplanes of $\sa$ intersecting $wC^A$ is 
    
    $$\H_w=\{x_{w(i)}-x_{w(j)}=1 | 1\leq i<j\leq n, w(i)<w(j)\}.$$
    
    There are two natural orders on the hyperplanes in $\H_w$; namely, hyperplane $H_1$ is less than hyperplane $H_2$ if all the points in $wC^A$ which are on the same  (opposite) side of $H_1$ as the origin are also on the same (opposite) side of $H_2$ as the origin. Thus, $\H_w$ can be considered as a poset. The  set of ceilings of some region of $\sa$ in $wC^A$ is an antichain of this poset. As Theorem \ref{bija} states below, the reverse is also true, and so the antichains of  $\H_w$  are in bijection with the regions of $\sa$ in $wC^A$. To avoid any confusion we now (re-)define the poset we consider.
    
    Let $$Q_w=\{(i, j) : 1\leq i<j\leq n, w(i)<w(j)\}$$  partially ordered by $$(i, j)\leq (r, s) \text{ if } r\leq i<j\leq s.$$   We think of $(i,j) \in Q_w$ as the hyperplane $x_{w(i)}-x_{w(j)}=1$ in $\H_w$.  Note that in $Q_w$  we have  $(i, j)\leq (r, s)$ if and only if 
    all  points in $wC^A$ which are on the same  side of $x_{w(i)}-x_{w(j)}=1$ as the origin are also on the same side of $x_{w(r)}-x_{w(s)}=1$ as the origin.
    
\begin{figure}[htbp] 
\begin{center} 
\includegraphics[scale=.5]{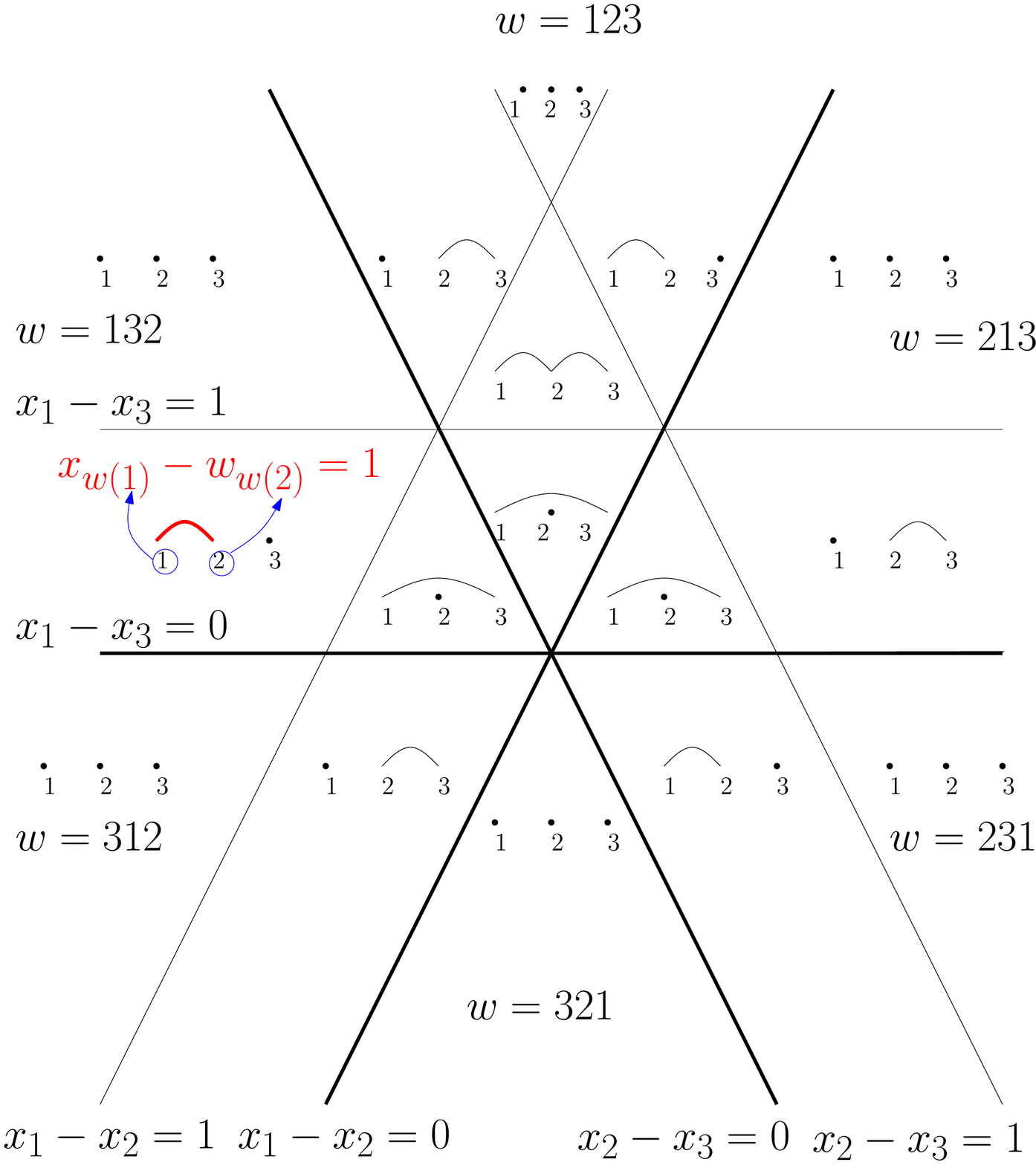} 
\caption{Labeling the regions of $\mathcal{S}_3^A$ by partitions.} 
\label{wnumomega} 
\end{center} 
\end{figure}

    We represent antichains of $Q_w$ as partitions of $[n]$, where we draw an arc $(i, j)$ in the diagram  if $(i, j) \in Q_w$ is in the antichain. For basic definitions about partitions  see \cite[Chapter 1]{ec1}. Bijecting the regions of $\sa$ to their set of ceilings, and the set of ceilings to the corresponding antichains in $Q_w$, which we represent as partitions,  we obtain a labeling of the regions of $\sa$ by partitions, as shown on Figure \ref{wnumomega}.

\begin{theorem} \label{bija}
\cite[Section 5]{sta0}, \cite[Theorem 2.1]{sta} 
  There is a bijection between the regions of $\sa$ contained in the cone $wC^A$ and the antichains of the poset $Q_w$. In particular,  \begin{equation} \label{a} |R(\sa)|=\sum_{w \in \S_n} j(Q_w), \nonumber \end{equation}
where    $j(Q_w)$ denotes the number of antichains of the poset $Q_w$.   \end{theorem} 
     
     \proof It is clear from the above that there is an injective map from  the regions of  $\sa$ to the multiset of the antichains of the posets $Q_w$, 
    $w \in \S_n$.  Since it is known that $|R(\sa)|=(n+1)^{n-1}$ \cite{shi} and  $\sum_{w \in \S_n} j(Q_w)=(n+1)^{n-1}$ can be proved without reference to $\sa$ (see Corollary \ref{cor:aa} in Section \ref{sec:aa}) the map also has to be surjective and  Theorem \ref{bija} follows.
     \qed
     
     \medskip

       Studying the relations of the poset $Q_w$ we see that the antichains of $Q_w$ correspond to nonnesting  $A_{n-1}$-partitions  if we think of $(k, l) \in Q_w$ as an arc in a partition of $[n]$.

   \subsection{Posets and the regions of $\scc$.} \label{sub:cc} Pick a region $R$ of $\scc$ in the cone $wC^{BC}$ of $\cb$, $w \in \S_n^B$. The set of hyperplanes of $\scc$ that  intersect $wC^{BC}$ is
   
   $$\H_w^C=\H_w^+\cup \H_w^-\cup \H_w^s,$$
  
  \noindent where 
  
  $$\H_w^-= \{ x_{w(i)}-x_{w(j)}=1 \mid i<j, 0<w(i)<w(j)\},$$ 
  
  $$\H_w^+= \{x_{w(i)}-x_{w(j)}=1 \mid i<j, w(j)<0<w(i)\},$$
  
  and $$\H_w^c= \{x_{w(i)}-x_{w(-i)}=1 \mid i<-i, w(-i)<0<w(i)\}.$$

 Taking into consideration that $w(i)=-w(-i)$ for all $i \in [\pm n]$, we can write

     $$\H_w^C= \{ x_{w(i)}-x_{w(j)}=1 \mid i<j, 0<w(i)\leq|w(j)|\},$$

   Note that if $x_{w(i)}-x_{w(j)}=1$, $i<j$, and $x_{w(i')}-x_{w(j')}=1$, $i'<j'$, belongs to $\H_w^C$  and $R$ is on the same side of the hyperplane $x_{w(i)}-x_{w(j)}=1$  as the origin and $i'\leq i<j\leq j',$ then $R$ is also on the same side of the hyperplane $x_{w(i')}-x_{w(j')}=1$ as the origin, since $x_{w(i')}-x_{w(j')}\leq x_{w(i)}-x_{w(j)}<1$. Considering all such implications among the hyperplanes of $\H_w$ we arrive to a partial order (there are two choices of partial order, pick one) on the hyperplanes. Define the poset  $$Q^C_w=\{(i, j), (-j, -i) \mid   i<j, 0<w(i)\leq|w(j)| \}$$ with the partial ordering inherited from the hyperplanes: $$(i, j)\leq (r, s) \text{ if } r\leq i<j\leq s.$$ We can think of mapping a region to its ceilings, or to its floors. In either case Theorem \ref{bijc} follows. For a related  bijection  between the positive chambers of the Shi arrangement and order ideals of the root poset of  corresponding type  see \cite[Theorem 5.1.13]{arm} and \cite{cel-pap}.

      \begin{theorem} \label{bijc}    The regions of $\scc$ contained in $wC^{BC}$ are in bijection with the antichains of ~$Q^C_w$.
   In particular,  \begin{equation} \label{a} |R(\scc)|=\sum_{w \in \S^B_n} j(Q^C_w), \nonumber \end{equation}
where    $j(Q^C_w)$ denotes the number of antichains of the poset $Q^C_w$.   \end{theorem} 
     
     \proof It is clear from the above that there is an injective map from  the regions of  $\scc$ to the multiset of the antichains of the posets $Q^C_w$, 
    $w \in \S^B_n$.  Since it is known that $| R(\scc)|={(2n+1)}^{n}$ \cite{shi2} and  $\sum_{w \in \S^B_n} j(Q^C_w)=(2n+1)^{n}$ can be proved without reference to $\scc$ (see Corollary \ref{cor:bb} in Section \ref{sec:bb})   the map also has to be surjective and  Theorem \ref{bijc} follows.
     \qed
     
\medskip

\begin{figure}[htbp] 
\begin{center} 
\includegraphics[scale=.85]{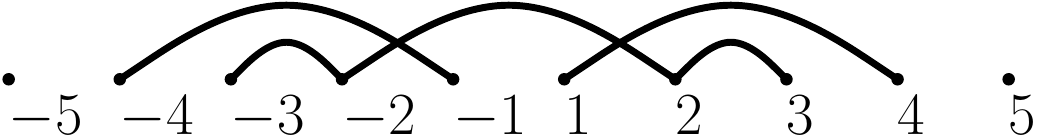} 
\caption{A $C_5$-partition.} 
\label{C-part} 
\end{center} 
\end{figure}

Note that the antichains of $Q^C_w$ correspond to nonnesting $C_n$-partitions   if we think of $(k, l) \in Q^C_w$ as an arc in a partition of $[\pm n]$. Recall that a \textbf{nonnesting  $C_n$-partition}  of $[\pm n]$ can be thought of as a nonnesting diagram of arcs, which are drawn over the ground set $-n, -n+1, \ldots, -2, -1, 1, 2, \ldots, n-1, n$ (in this order) such that if there is an arc between $i$ and $j$, for $i, j, \in [\pm n]$, then there is also an arc between $-j$ and $-i$ (there are no multiple arcs). See Figure \ref{C-part} for an example.

  Bijecting the regions of $\scc$ to their set of ceilings, and the set of ceilings to the corresponding antichains in $Q^C_w$, which we represent as $C_n$-partitions,  we obtain a labeling of the regions of $\scc$ by partitions, analogously to the type $A_{n-1}$ case.

  \section{Sequences and Shi arrangements in type $A_{n-1}$}
    \label{secaa}

 In this section we construct a bijection between the regions of $\sa$ and $(n+1)$-tuples of sequences $a_1\ldots a_n$, $a_i \in [n+1]$, for $i \in [n]$, such that every such  sequence appears in exactly one of the  $(n+1)$-tuples.  Furthermore, exactly one among the $n+1$  sequences assigned to a region is a parking function, thereby also leading to a bijection between the regions of $\sa$ and parking functions. The same  bijection  previously appeared in the paper by Athanasiadis and Linusson \cite{ath-lin}. Our exposition makes  the enumeration of the ceiling and floor statistic on the regions on $\sa$ transparent, and that it readily generalizes to bijections in the type $C_n$ case. The ceiling and floor statistics on Shi arrangements was also used and studied by Armstrong and Rhoades in their beautiful paper on the Shi and Ish arrangements \cite{arm-rho}.  The  ideas of this section appear explicitly or implicitly in \cite{ath-lin} and  \cite{arm-rho}.

  For ease of exposition we consider $n+1$ copies of the arrangement $\sa$, denoted by $${(\sa)}^{\sqcup (n+1)}:={(\sa)}^{(1)}\sqcup \cdots \sqcup {(\sa)}^{(n+1)},$$ and biject the regions of  ${(\sa)}^{\sqcup (n+1)}$ defined as the  regions of  ${(\sa)}^{(1)}, \ldots, {(\sa)}^{(n)}$ and ${(\sa)}^{(n+1)}$ with the sequences $a_1\ldots a_n$, $a_i \in [n+1]$, for $i \in [n]$.

 The \textbf{type} of an $A_{n-1}$-partition $\pi$ is the integer partition $\l$ whose parts are the sizes 
  of the blocks of $\pi$. 

\begin{theorem} \label{a1} \cite{ath} There is a bijection $b$ between the set of type $\l=(\l_1, \ldots, \l_d)$ nonnesting $A_{n-1}$-partitions and pairs $(S, g)$, where $S$ is a $d$-subset of $[n]$ and the map $g:S\rightarrow \{\l_1, \ldots, \l_d\}$ is such that $|g^{-1}(i)|=r_i$, $0 \leq i$. 

\end{theorem}

\proof  Given a  type $\l=(\l_1, \ldots, \l_d)$ nonnesting $A_{n-1}$-partition, let $S$ be the leftmost elements of its blocks. Let $g(s)=k$ if  $s \in S$ is in a block of size $k$. It can be shown by induction on $n$ that  the set $S$ and function $g$ defined this way uniquely determine the nonnesting partition they came from. \qed

\medskip

\begin{figure}[htbp] 
\begin{center} 
\includegraphics[scale=.2]{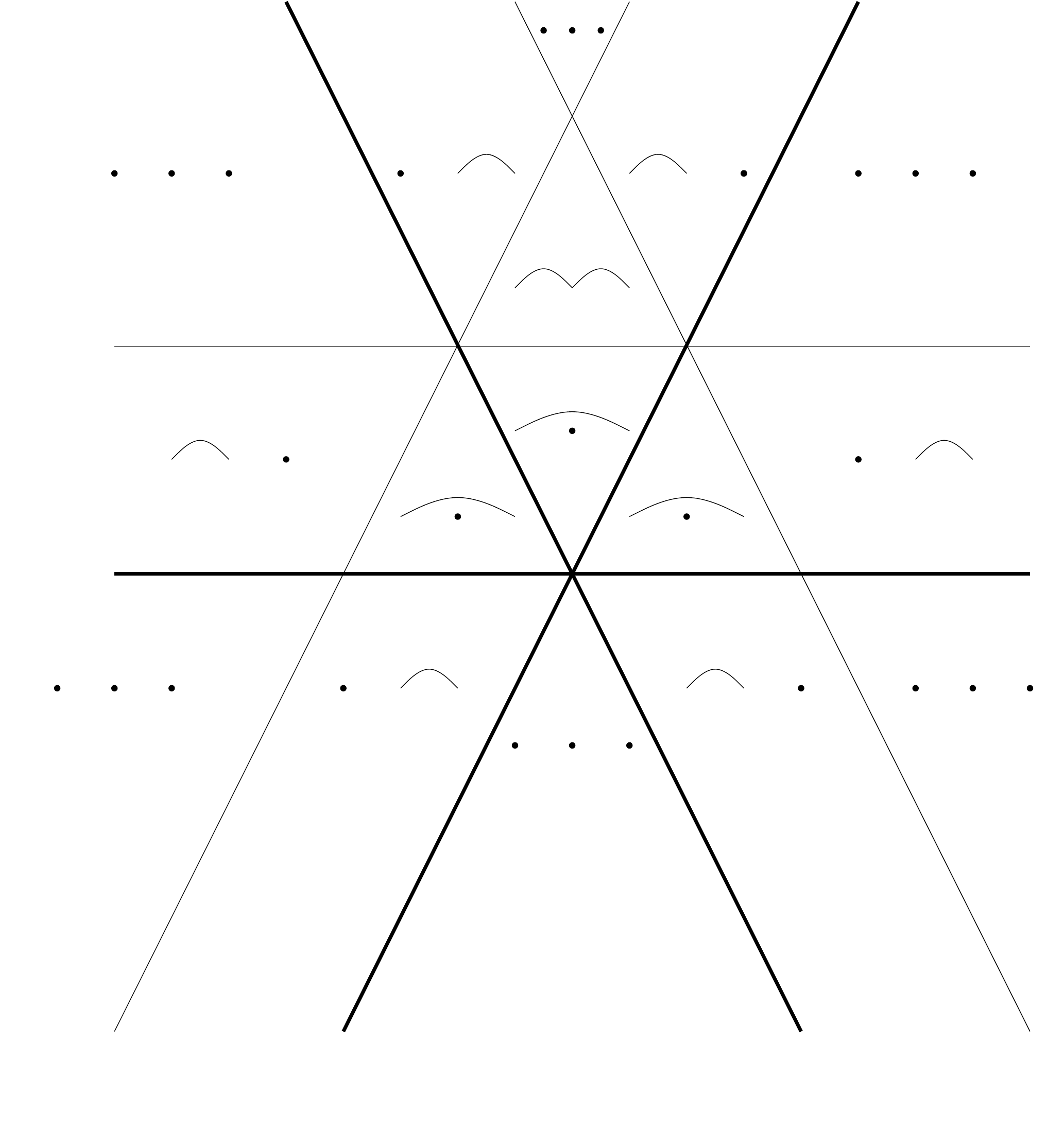} 
\includegraphics[scale=.2]{n+1.pdf} 
\includegraphics[scale=.2]{n+1.pdf} 
\includegraphics[scale=.2]{n+1.pdf} 
\caption{Labeling of the regions of $(\mathcal{S}_3^A)^{\sqcup 4}$ by nonnesting partitions.} 
\label{"n+1"} 
\end{center} 
\end{figure}

\begin{figure}[htbp] 
\begin{center} 
\includegraphics[scale=.5]{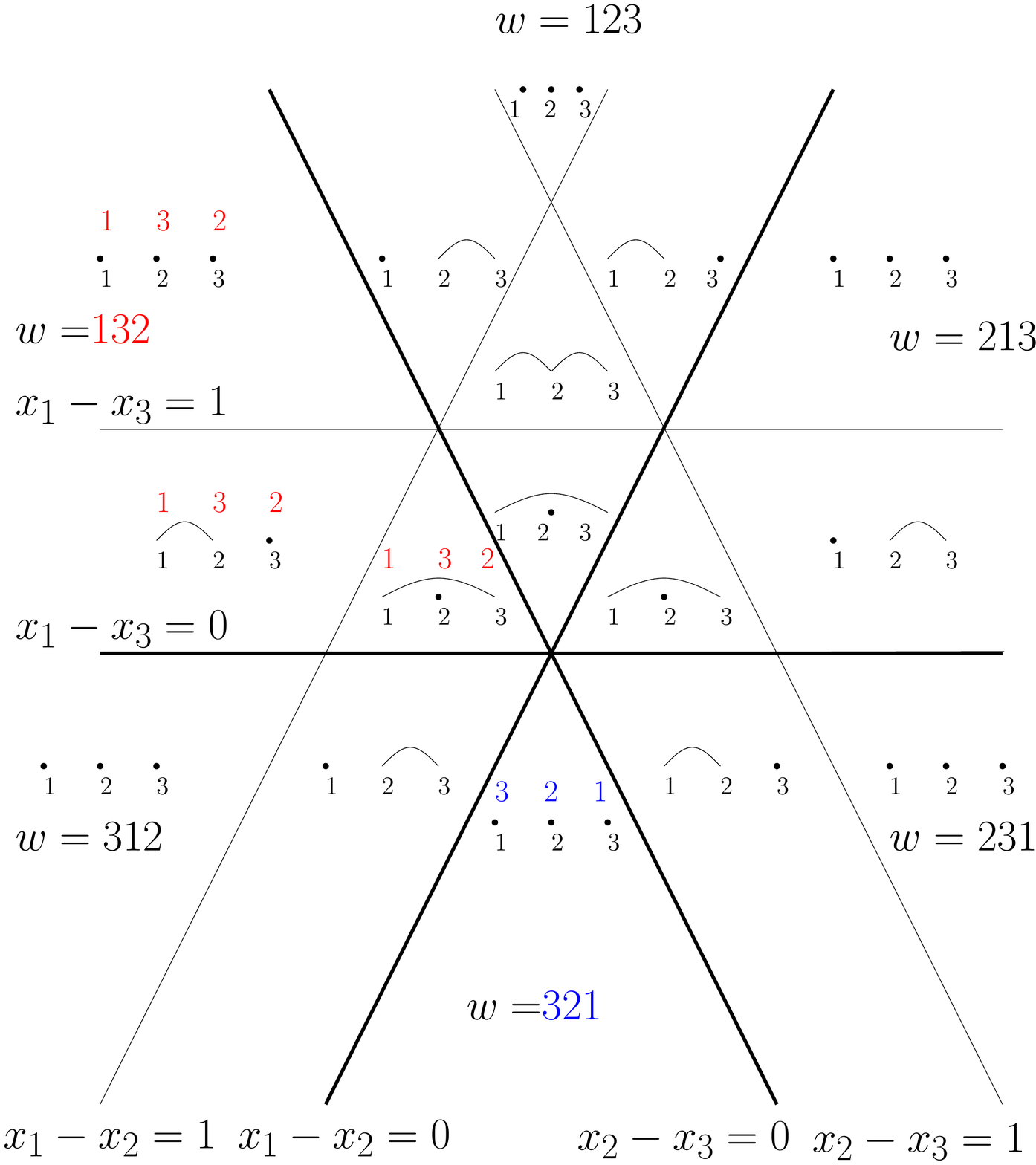} 
\caption{The $w$-labels of partitions labeling of the regions of $\mathcal{S}_3^A$ are shown for the permutations $132$ and $321$.} 
\label{wlabel} 
\end{center} 
\end{figure}

Label each region of $\san$ by the nonnesting $A_{n-1}$-partition corresponding to an antichain of $Q_w$, $w \in \S_n$, as described in Section \ref{sub:aa} and shown on Figure \ref{"n+1"}. Each region of  $\san$ is completely specified by a number $k \in [n+1]$ (specifying which copy of $\sa$  we are in in $\san$),  a permutation $w \in \S_n$ (specifying the cone of $\sa$), and a nonnesting $A_{n-1}$-partition $\pi$ (specifying the ceilings of $R$ in $wC^A$). While we generally think of $\pi$ as on the vertices $1, 2, \ldots, n-1, n$, in this order, the $A_{n-1}$-partition $\pi$ also has \textbf{$w$-labels} $ w(1), \ldots, w(n-1), w(n)$. See Figure \ref{wlabel}.

\begin{lemma} \label{a2}The number of regions of $\sa$ containing the nonnesting $A_{n-1}$-partition $\pi$ of type $\l$ is equal to 

\begin{equation} \label{acc} {n \choose {\l_1, \ldots, \l_d}}. \end{equation}
\end{lemma}

\proof In this proof we effectively count the number of  permutations $w\in \S_n$ such that $\pi$ is an antichain in the poset $Q_w$, since the latter is equal to the of number regions of $\sa$ containing the nonnesting $A_{n-1}$-partition $\pi$. Given a nonnesting $A_{n-1}$-partition $\pi$ of type $\l$ there are   ${n \choose {\l_1, \ldots, \l_d}}$ ways to choose  the values of the $w$-labels which go into the blocks of $\pi$. Since in each  block the $w$-labels  increase,  equation \eqref{acc} follows. 
\qed

\medskip

Given a type $\l=(\l_1, \ldots, \l_d)$ nonnesting $A_{n-1}$-partition  $\pi$, denote by $S_{\pi}$ the set and $g_{\pi}$ the function   from Theorem \ref{a1}. Let $\bar{S}(\pi)$ be the multiset consisting of $\l_i$  copies of each element of $g_{\pi}^{-1}(\l_i)$, for each part in the set (not multiset!) $\{\l_1, \ldots, \l_d\}$. An \textbf{$n$-shifted permutation} of $\bar{S}(\pi)$ is a permutation of the elements of the multiset $\bar{S}(\pi)$ such that each  entry is increased by $k \in \{0, 1, \ldots, n\}$ and taken modulo $n+1$. For example the $2$-shifted permutations of $\{\{1, 2\}\}$ are $12, 21, 23, 32, 31, 13$.

\begin{theorem} \label{a3} There is a bijection $\phi$ between the regions of $\san$ labeled by the nonnesting $A_{n-1}$-partition $\pi$ of type $\l$ and $n$-shifted permutations of the multiset $\bar{S}(\pi)$. 

\end{theorem}

\proof There are multiple ways to set up this bijection. We present one way here and note how to define a family of bijections satisfying Theorem \ref{a3}. 

Given  the nonnesting $A_{n-1}$-partition $\pi$, a  permutation $w$ for which $\pi$ is an antichain in $Q_w$, and an integer  $k \in [n+1]$ specifying which copy of $\sa$  we are in in $\san$, order the blocks  of $\pi$  by increasing size. The blocks of the same size are ordered lexicographically according to the $w$-labels on them. Order the numbers in the multiset  $\bar{S}(\pi)$ so that the numbers with less multiplicites are smaller. Among the numbers with the same multiplicity order them according to the natural order on integers. The previous two orders yield a correspondence $b_3$ between the blocks $B$ of $\pi$ and the numbers from $\bar{S}(\pi)$.  (This correspondence could of course  be defined in several ways leading to different bijections.) Let the values of the $w$-labels of $B$ specify the positions that the number $b_3(B)+k-1 \textrm{ mod } n+1$ is taking.

Correspondence $b_3$ could also be naturally defined by the bijection given in Theorem \ref{a1}. As it turns out both descriptions of $b_3$  in the type $A_{n-1}$ case are the same.

\begin{figure}[htbp] 
\begin{center} 
\includegraphics[scale=.5]{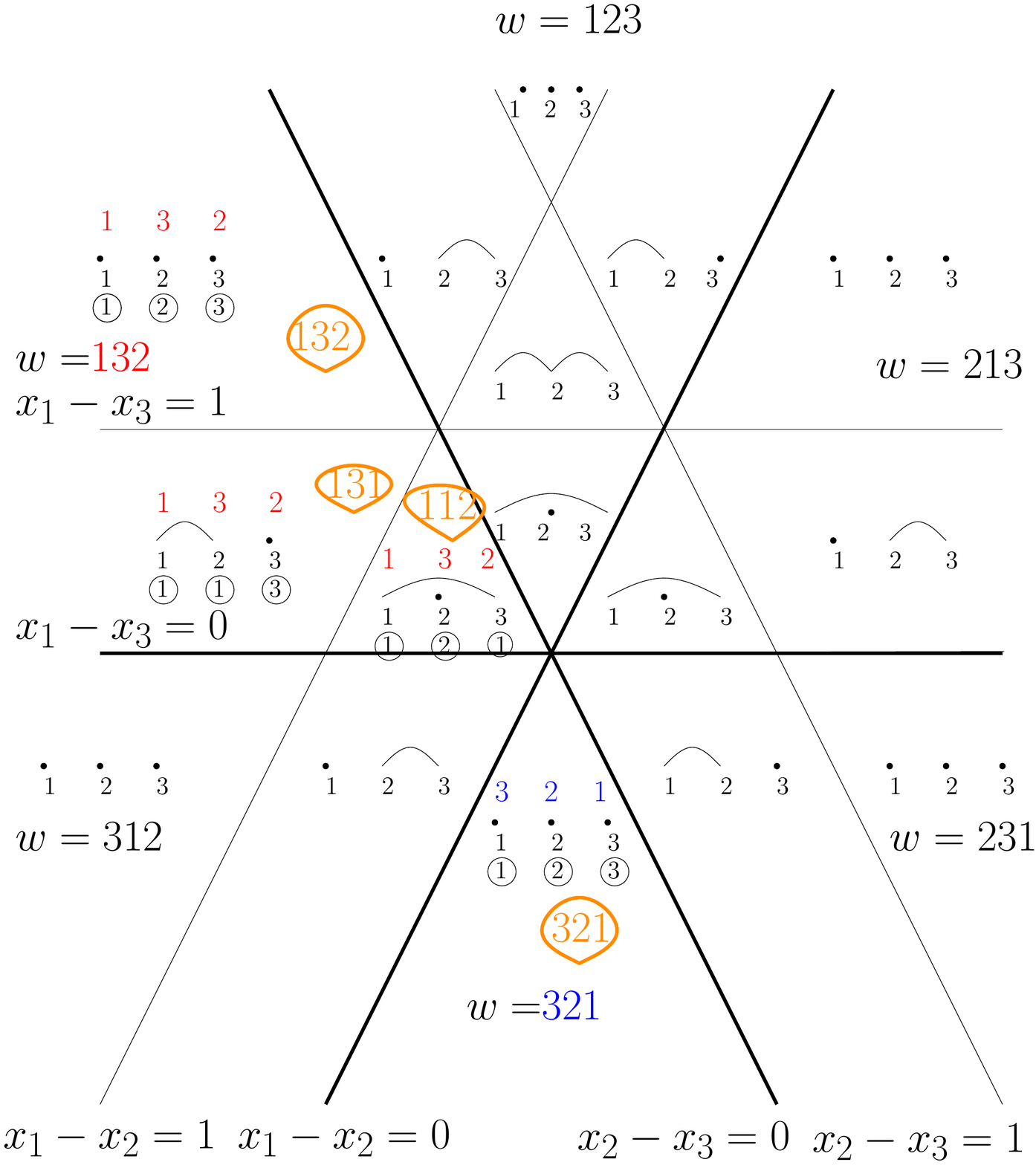} 
\caption{Constructing the bijection on $(\mathcal{S}^A_3)^{(1)}$.} 
\label{fig:bij} 
\end{center} 
\end{figure}

Figure \ref{fig:bij} shows the construction of the bijection on $(\mathcal{S}^A_3)^{(1)}$. For example, consider the region $R$ in $(132)C^A$ with ceiling $x_1-x_3=x_{w(1)}-x_{w(2)}=1$.  The $w$-labels are written above the partition and the numbers corresponding to the elements of the blocks are below the partition and are circled individually. Then, to get the sequence corresponding to the region, read the circled numbers in the order specified by the $w$-labels. The resulting $3$-digit sequence is circled on  Figure \ref{fig:bij}.

\begin{figure}[htbp] 
\begin{center} 
\includegraphics[scale=.5]{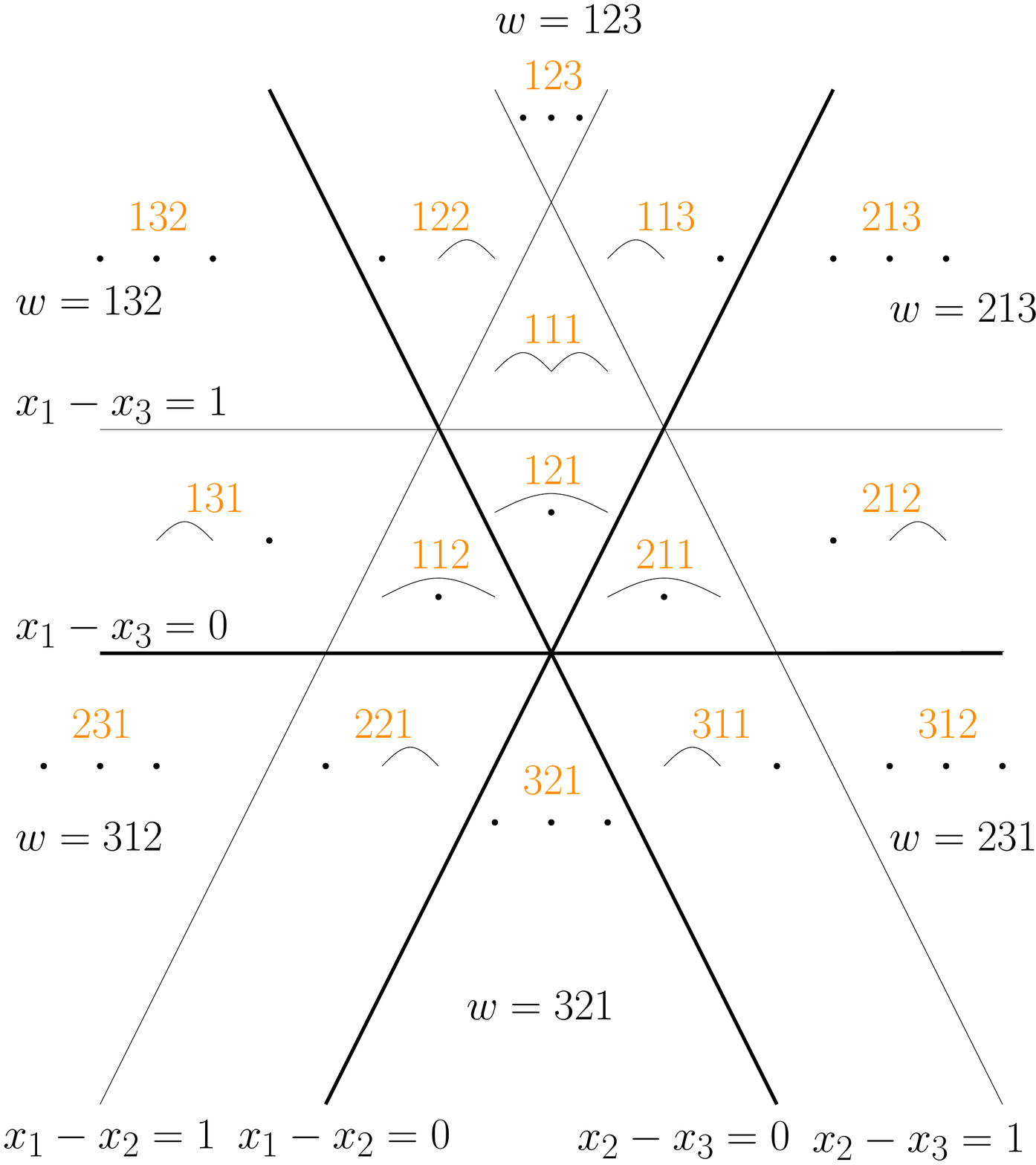} 
\caption{The bijection on $(\mathcal{S}^A_3)^{(1)}$.} 
\label{bijcomplete} 
\end{center} 
\end{figure}

For the restriction of the bijection to $(\mathcal{S}^A_3)^{(1)}$ see Figure \ref{bijcomplete}.

The above defined map is a  bijection  between the regions of $\san$ labeled by the nonnesting $A_{n-1}$-partition $\pi$ of type $\l$ and  $n$-shifted permutations of the multiset $\bar{S}(\pi)$, which can be shown by writing down an explicit inverse, or by noting that it is injective and the domain and codomain are equinumerous.
\qed

\medskip

Extend the map $\phi$ defined in the proof of  Theorem \ref{a3} to a map between 
all regions of $\san$ and the set of sequences $\A(n)=\{a_1\ldots a_n | a_i \in [n+1], i \in [n]\}$.

\begin{theorem} \label{a4} (cf. \cite{ath-lin})  The map  $\phi: R(\san)\rightarrow \A(n)$ is a bijection.  
\end{theorem}

\begin{theorem} \label{a5} (cf. \cite{ath-lin})  The restriction of the bijection $\phi$ to the first copy of the Shi arrangement ${\sa}^{(1)}$ is a bijection between the regions of the Shi arrangement and parking functions.\end{theorem} 

We leave the details of the proofs of Theorems \ref{a4} and \ref{a5} to the reader. Hint: see \cite[Exercise 5.49]{ec2}. 

\begin{corollary} 
\begin{equation} \label{uhhu}
  \sum_{\l \vdash n} {n+1 \choose d} \frac{d!}{m_{\l}}{n \choose {\l_1, \l_2, \ldots, \l_d}}=(n+1)^n,\end{equation} 
  where $m_{\l}=\prod_{i=1}^n r_i!$, if $r_i$ denotes the number of parts of $\l$ equal to $i$.

   \end{corollary}

\proof Kreweras 
\cite[Theorem 4]{kre}  proved that the number of noncrossing partitions of $[n]$ of type $\l$ is equal  $$ \frac{n!}{m_{\l}(n-d+1)!},$$ 
 
 \noindent where $d$ denotes the number of parts of $\l$. Athanasiadis \cite[Theorem 3.1]{ath} gave a bijection between noncrossing and nonnesting $A_{n-1}$-partitions which preserves type. Thus, the total number of nonnesting $A_{n-1}$-partitions of type $\l$ labeling the regions of $\san$ is 
 \begin{equation}\label{n+1} (n+1)\frac{n!}{m_{\l}(n-d+1)!}= {n+1 \choose d} \frac{d!}{m_{\l}}.\end{equation}
 Equation \eqref{n+1} together with Lemma \ref{a2} and Theorem \ref{a4} imply equation \eqref{uhhu}. 
\qed

\medskip

Theorem \ref{gf} from the introduction is  a corollary of the proofs of Theorems \ref{a3}, \ref{a4} and \ref{a5}.  For further details see Section \ref{sec:aa}, and in particular Theorem \ref{bij}.
\medskip

\noindent \textbf{Theorem 2.}
$$ \sum_{R \in R(\sa)} q^{c(R)}= \sum_{R \in R(\sa)} q^{f(R)}=\frac{1}{n+1}\sum_{\textrm{a} \in \A(n)} q^{n-d(\a)}=\sum_{\a \in PF(n)} q^{n-d(\a)}.$$

    \section{Sequences and Shi arrangements in type $C_{n}$}
    \label{secc}

 In this section we construct a  bijection  between the regions of  $\scc$ and the set of  sequences $\A^{C}(n)=\{a_1\ldots a_n|  a_i \in  [\pm n] \cup \{0\}, i \in [n]\}$.   Our proof  yields  enumeration of regions by  the ceiling and floor statistic, which we express in a generating function form. 
  
\begin{figure}[htbp] 
\begin{center} 
\includegraphics[scale=.85]{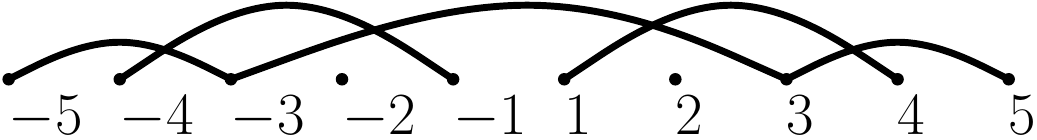} 
\caption{A type $(2, 1)$ nonnesting $C_5$-partition.} 
\label{C-nonnest} 
\end{center} 
\end{figure}

 The \textbf{type} of a $C_{n}$-partition $\pi$ is the integer partition $\l$ whose parts are the sizes 
  of the nonzero blocks of $\pi$, including one part for each pair of blocks $\{B, -B\}$.   The  zero block is a block $B$ such that $B=-B$. Figure \ref{C-nonnest} shows a nonnesting $C_5$-partition with blocks $\{2\}, \{-2\}, \{-1, -4\}, \{1, 4\}, \{-5,-3,$ $ 3, 5\}$.  The last block is a zero block, and so the type of this partition is $(2, 1)$.

 The following theorem is based on  a bijection of Fink and Iriarte \cite{fin-iri}  between noncrossing and nonnesting $C_n$-partitions which preserves type and a bijection of Athanasiadis \cite{ath} between noncrossing $C_n$-partitions and pairs $(S, g)$, where $S$ is a set and $g$ is a function subject to the conditions stated below.

\begin{theorem} \label{1}There is a bijection between the set of type $\l=(\l_1, \ldots, \l_d)$ nonnesting $C_n$-partitions and pairs $(S, g)$, where $S$ is a $d$-subset of $[n]$ and the map $g:S\rightarrow \{\l_1, \ldots, \l_d\}$ is such that $|g^{-1}(i)|=r_i$, $0 \leq i$. 

\end{theorem}

\proof \cite[Theorem 2.4]{fin-iri} establishes a type-preserving bijection $b_1$ between nonnesting and noncrossing $C_n$-partitions, and \cite[Theorem 2.3]{ath} provides a bijection $b_2$ between the set of type $\l=(\l_1, \ldots, \l_d)$ noncrossing $C_n$-partitions and pairs $(S, g)$, where $S$ is a $d$-subset of $[n]$ and the map $g:S\rightarrow \{\l_1, \ldots, \l_d\}$ is such that $|g^{-1}(i)|=r_i$, $0 \leq i$. \qed

\medskip

Label each region of $\scc$ by the nonnesting $C_n$-partition corresponding to an antichain of $Q_w^C$, $w \in \S_n^B$, as described in Section \ref{sub:cc}. Each region of  $\scc$ is completely specified by a nonnesting $C_n$-partition $\pi$ and $w \in \S_n^B$. While we generally think of $\pi$ as on the vertices $-n, -n+1, \ldots, -1, 1, 2, \ldots, n-1, n$, in this order, the $C_n$-partition $\pi$ also has \textbf{$w$-labels} $w(-n), w(-n+1), \ldots, w(-1), w(1), \ldots, w(n-1), w(n)$. 

\begin{lemma} \label{2}The number of regions of $\scc$ containing the nonnesting $C_n$-partition $\pi$ of type $\l$ is equal to 
\begin{equation} \label{cc} {n \choose {\l_1, \ldots, \l_d, n-|\l|}} \prod_{i=1}^d 2^{\l_i}. \end{equation}
\end{lemma}

\proof In this proof we effectively count the number of signed permutations $w\in \S_n^B$ such that $\pi$ is an antichain in the poset $Q_w^C$, since the latter is equal to the number regions of $\scc$ containing the nonnesting $C_n$-partition $\pi$. Given a nonnesting $C_n$-partition $\pi$ of type $\l$ there are   ${n \choose {\l_1, \ldots, \l_d}, n-|\l|}$ ways to choose  the absolute values of the $w$-labels which go into the blocks of $\pi$. Let $2(n-\l)$ be the  size of the zero block of $\pi$. The signs and order of the $w$-labels in the zero block of $\pi$ are determined: there have to be $(n-\l)$ positive numbers in increasing order followed by their $(n-\l)$ negatives in increasing order. Each nonzero  block is comprised  of  a possibly empty sequence of positive $w$-labels in increasing order followed by a possibly empty sequence of negative $w$-labels in increasing order. There are exactly $2^{\l_i}$ ways to decide the signs among $\l_i$ numbers, and once the signs are decided so is the order. Thus, equation \eqref{cc} follows. 
\qed

\medskip

Given a type $\l=(\l_1, \ldots, \l_d)$ nonnesting $C_n$-partition  $\pi$, denote by $S_{\pi}$ the set and $g_{\pi}$ the function   from Theorem \ref{1}. Let $\bar{S}(\pi)$ be the multiset consisting of $n-\l$ $0$'s, and $\l_i$  copies of each element of $g_{\pi}^{-1}(\l_i)$, for each part in the set (not multiset!) $\{\l_1, \ldots, \l_d\}$. A \textbf{marked permutation} of $\bar{S}(\pi)$ is a permutation of the elements of the multiset $\bar{S}(\pi)$ such that each  nonzero  entry has a $\pm$ sign in addition. For example the marked permutations of $\{\{0, 1, 1\}\}$ are $011, 101, 110, 0-1-1, -10-1, -1-10, 01-1, 10-1, 1-10, 0-11, -101, -110$ (we omitted the $+$ signs).

\begin{theorem} \label{3} There is a bijection $\phi$ between the regions of $\scc$ labeled by the nonnesting $C_n$-partition $\pi$ of type $\l$ and marked permutations of the multiset $\bar{S}(\pi)$. 

\end{theorem}

\proof There are multiple ways to set up this bijection. We present two natural  ways here and note how to define a family of bijections satisfying Theorem \ref{3}. 

Given  the nonnesting $C_n$-partition $\pi$ and a signed permutation $w$ for which $\pi$ is an antichain in $Q_w^C$, order the (pair of) blocks $\{B, -B\}$ of $\pi$ as follows. If there is a zero block, then it comes first. The other blocks are ordered by increasing size ($|B|$), and the blocks of the same size are ordered lexicographically according to the $w$-labels on them (on $B$s). Order the numbers in the multiset  $\bar{S}(\pi)$ so that the $0'$s come first, and among the other numbers the numbers with less multiplicites are smaller. Among the numbers with the same multiplicity order them according to the natural order on integers. The previous two orders yield a correspondence $b_3$ between the blocks $B$ of $\pi$ and the numbers from $\bar{S}(\pi)$.  (This correspondence could of course  be defined in several ways leading to different bijections.) Let the absolute values of the $w$-labels of $B$ specify the positions that the number $b_3(B)$ is taking. For the $w$-labels of nonzero blocks which are negative add a $-$ to the number in the corresponding spot.

Correspondence $b_3$ could also be naturally defined by the bijection given in Theorem \ref{1}. 

The above defined maps are  bijections  between the regions of $\scc$ labeled by the nonnesting $C_n$-partition $\pi$ of type $\l$ and marked permutations of the multiset $\bar{S}(\pi)$, which can be shown by writing down explicit inverses, or by noting that they are injective and the domains and codomains are equinumerous.
\qed

\medskip

Extend the map $\phi$ defined in the proof of  Theorem \ref{3} to a map between 
all regions of $\scc$ and the set of sequences $\A^{C}(n)=\{a_1\ldots a_n | a_i \in [\pm n]\cup \{0\}, i \in [n]\}$, to obtain the following corollaries as in the type $A_{n-1}$ case.

\begin{theorem} \label{4} The map  $\phi: R(\scc)\rightarrow \A^{C}(n)$ is a bijection.  
\end{theorem}

\begin{corollary} 
$$\sum_{\l \vdash n} \frac{n!}{m_{\l}(n-d)!} {n \choose {\l_1, \ldots, \l_d, n-|\l|}} \prod_{i=1}^d 2^{\l_i}=(2n+1)^n.$$
\end{corollary}

\proof
Athanasiadis \cite{ath} proved that the number of  nonnesting   $C_n$-partitions  of type $\l$ is  $$ \frac{n!}{m_{\l}(n-d)!},$$ which together with Lemma \ref{2} and Theorem \ref{4} imply the above equality. 
\qed

\medskip

Theorem \ref{bgf} is   a corollary of the proofs of  Theorems \ref{1}, \ref{3} and  \ref{4} (use the second definition of $b_3$ in the proof of Theorem \ref{3}). For further details see Section \ref{sec:bb}, and in particular Theorem \ref{cbij}.

\medskip

\noindent \textbf{Theorem 4.}
$$ \sum_{R \in R(\scc)} q^{c(R)}= \sum_{R \in R(\scc)} q^{f(R)}=\sum_{\a \in \A^{C}(n)} q^{n-d^{C}(\a)}.$$

 \section{Posets and sequences in  type $A_{n-1}$}
    \label{sec:aa}

In this section we revisit the type $A_{n-1}$ world of  posets $Q_w$, $w \in \S_n$, and parking functions of length $n$ and state their relation explicitly without the mention of arrangements. Much of the considerations of this section appear in the work of Athanasiadis and Linusson \cite{ath-lin} and Armstrong and Rhoades \cite{arm-rho} either explicitly or implicitly. We highlight our perspective on the relation of the posets and sequences and study their properties in detail.  We carry out a similar agenda for the posets $Q^{C}_w$, $w \in \S_n^B$, and sequences in $\A^{C}(n)$ in the next section.

 Recall that $$Q_w=\{(i, j) : 1\leq i<j\leq n, w(i)<w(j)\}$$  is partially ordered by $$(i, j)\leq (r, s) \text{ if }r\leq i<j\leq s.$$ We explore the refinements of the equation
 
 \begin{equation} \label{q} \sum_{w \in \S_n} j(Q_w)=(n+1)^{n-1},\end{equation}

\noindent which follows from  Theorems \ref{bija} and \ref{a4}.  In the process we reiterate the  proof of equation \eqref{q} without reference to arrangements. 

 Partition the set of parking functions of length $n$, $PF(n)$, according to the cardnality of the set $\{a_1, a_2, \ldots, a_n\}$. Let $S_k(n)=\{(a_1, a_2, \ldots, a_n) \in PF(n) : |\{a_1, a_2, \ldots, a_n\}|=k\}$. Then
 
 $$PF(n)=\bigcup_{k=1}^n S_k(n).$$

 Partition the multiset of antichains $\mathcal{M}(n)$ of $Q_w$, $w \in \S_n$, according to the cardinality of the antichains.   Let $M_k(n)=\{\{ \{(i_1, j_1), \ldots, (i_k, j_k)\} \in \mathcal{M}(n) \}\}$. Then
 
 $$\mathcal{M}(n)=\bigcup_{k=0}^{n-1} M_k(n).$$

The following theorem can be deduced  from the work of Athanasiadis and Linusson \cite{ath-lin} and Armstrong and Rhoades \cite{arm-rho}. 

\begin{theorem} \label{bij}
 
 $$|S_k(n)|=|M_{n-k}(n)|, \mbox{ } k \in [n].$$ 
 \end{theorem}
 
We prove Theorem \ref{bij} by providing a bijection between the sets $S_k(n)$ and $M_{n-k}(n)$, $k \in [n]$. Before proceeding to the proof of  Theorem \ref{bij}   we partition the sets  $S_k(n)$ and $M_{n-k}(n)$, $k \in [n]$, further. 
 
 Partition the set of parking functions of length $n$ with $k$ distinct numbers $S_k(n)$, $k \in [n]$, according to the $k$ distinct numbers appearing in the  sequence $(a_1, a_2, \ldots, a_n)$, and the number of times they appear. If $\{a_1, a_2, \ldots, a_n\}=\{c_1<c_2<\cdots<c_k\}$ and $c_i$ appears $o_i$ times in $(a_1, a_2, \ldots, a_n)$, $ i \in [k]$, let $$S_k^{\c, \o}(n)=\{(a_1, a_2, \ldots, a_n) \in S_k(n) |  
 \{\{a_1, a_2, \ldots, a_n\}\}=\cup_{i=1}^{k}\cup_{j=1}^{o_i}\{\{c_i\}\}  \mbox{ } \},$$
 
 \noindent where $\c=(c_1< \ldots< c_k)$, $\o=(o_1, \ldots, o_k)$, $o_i>0$, for $i\in [k]$, and $\sum_{i=1}^k o_i=n$.
 
 Given an antichain $\aa=\{(i_1, j_1), \ldots, (i_{n-k}, j_{n-k})\} \in M_{n-k}(n)$, $k \in [n]$, it naturally corresponds to  a nonnesting partition $\pi_{\aa}$  of $[n]$ with $k$ blocks,  where the arc diagram of $\pi_{\aa}$ consists of the arcs $(i_1, j_1), \ldots, (i_{n-k}, j_{n-k})$. Order the $k$ blocks of $\pi_{\aa}$ according to their smallest elements $c_i$, $i \in [k]$,  $c_1< \ldots< c_k$. Let $o_i>0$, $i \in  [k]$, be the number of elements in the $i^{th}$ block of $\pi_{\aa}$. Denote by $c(\aa)=(c_1< \ldots< c_k)$ and  $o(\aa)=(o_1, \ldots, o_k)$.
 Partition the multiset of antichains of length $n-k$ of  $Q_w$, $w \in \S_n$, $M_{n-k}(n)$, $k \in [n]$, according to $\c=(c_1< \ldots< c_k)$, $\o=(o_1, \ldots, o_k)$, $o_i>0$, for $i\in [k]$, and $\sum_{i=1}^k o_i=n$, as described above. Let 

$$M_{n-k}^{\c, \o}(n)=\{\{\aa=\{(i_1, j_1), \ldots, (i_{n-k}, j_{n-k})\} \in M_{n-k}(n)|  
c(\aa)=\c, o(\aa)=\o \}\},$$
 
 \noindent where $\c=(c_1< \ldots< c_k)$, $\o=(o_1, \ldots, o_k)$, $o_i>0$, for $i\in [k]$, and $\sum_{i=1}^k o_i=n$.

\begin{lemma} \label{un} \cite{ath}
The vectors $c(\aa)=\c$ and $o(\aa)=\o$, where  $\c=(c_1< \ldots< c_k)$, $\o=(o_1, \ldots, o_k)$, $k \in [n]$, $o_i>0$, for $i\in [k]$, $\sum_{i=1}^k o_i=n$, $c_1=1$, and $c_i \in \{c_{i-1}+1, c_{i-1}+2, \ldots, o_1+\cdots+o_{i-1}+1\}$, for $i\in \{2, \ldots, k\}$,  uniquely determine the antichain $\aa$.
\end{lemma}

 The following theorem can be deduced from the work of Athanasiadis and Linusson \cite{ath-lin} and Armstrong and Rhoades \cite{arm-rho}. 
 
 \begin{theorem} \label{co}
 $$|S_k^{\c, \o}(n)|=|M_{n-k}^{\c, \o}(n)|={ n \choose {o_1, \ldots, o_k}}, $$
 
 \noindent where $k \in [n]$, $\c=(c_1< \ldots< c_k)$, $\o=(o_1, \ldots, o_k)$, $o_i>0$, for $i\in [k]$, and $\sum_{i=1}^k o_i=n$.
 \end{theorem}
  
  \proof A bijective proof  can be given using Theorem \ref{a1} and the ideas of Theorem \ref{a3}. The enumeration is in Lemma \ref{a2}.  Note that arrangements do not enter the proof. \qed
    \medskip
  
 \noindent \textit{Proof of Theorem \ref{bij}.} Straightforward corollary of Theorem \ref{co}, since 
$$S_{k}(n)=\sum_{\c, \o} S_{k}^{\c, \o}(n)=\sum_{\c, \o} M_{n-k}^{\c, \o}(n)=M_{n-k}(n),$$  where $\c=(c_1< \ldots< c_k)$, $\o=(o_1, \ldots, o_k)$, $k \in [n]$, $o_i>0$, for $i\in [k]$, $\sum_{i=1}^k o_i=n$.
 \qed

 \begin{corollary} \label{cor:aa}
 
$$ \sum_{w \in \S_n} j(Q_w)=(n+1)^{n-1}.$$
 \end{corollary}
 
 \proof Theorems \ref{co} and  \ref{bij} extend to a bijection between $$MA(n)=\cup_{k=1}^n \cup_{\c, \o}M_{n-k}^{\c, \o}(n) \mbox{ and }PF(n)=\cup_{k=1}^n \cup_{\c, \o}S_{k}^{\c, \o}(n),$$ the cardinalities of which are  $ \sum_{w \in \S_n} j(Q_w)$ and $(n+1)^{n-1}$, respectively.
 \qed
 
 \section{Posets and sequences in type  $C_{n}$}
    \label{sec:bb}

In this section we revisit the type   $C_{n}$ world of  posets $Q^{C}_w$, $w \in \S_n^{B}$, and  sequences in $\A^{C}(n)$  and state their relation explicitly without the mention of arrangements.       
       
       Recall that 
       
$$Q^C_w=\{(i, j), (-j, -i) \mid   i<j, 0<w(i)\leq|w(j)| \}$$       
         is partially ordered by $$(i, j)\leq (r, s) \text{ if }r\leq i<j\leq s.$$ We explore the refinements of the equation
 
 \begin{equation} \label{qc} \sum_{w \in \S_n^B} j(Q^{C}_w)=(2n+1)^{n},\end{equation}

\noindent which follows from  Theorems \ref{bijc} and \ref{4}.  In the process we reiterate the  proof of equation \eqref{qc} without reference to arrangements. 

 Partition $\A^{C}(n)$  according to the number of nonzero absolute values in  the set $\{a_1, a_2, \ldots, a_n\},$ denoted by $d^{C}(\a)$ for $\a=a_1a_2\ldots a_n$. Let $S^{C}_k(n)=\{(a_1, a_2, \ldots, a_n) \in \A^{C}(n) : d^{C}(\a)=k\}$. Then
 
 $$\A^{C}(n)=\bigcup_{k=0}^n S^{C}_k(n).$$

 Partition the multiset of antichains $\mathcal{M}^{C}(n)$ of $Q^{C}_w$, $w \in \S^{B}_n$, according to the number of pairs $(i, j), (-j, -i)$, $i \leq j$,  in the antichains. Denote by $p(\pi)$ the number of pairs $(i, j), (-j, -i)$, $i \leq j$,  in the antichain $\pi \in \mathcal{M}^{C}(n)$.   Let $M^{C}_k(n)=\{\{ \pi \in \mathcal{M}(n) | p(\pi)=k \}\}$. Then
 
 $$\mathcal{M}^{C}(n)=\bigcup_{k=0}^{n} M^{C}_k(n).$$

\begin{theorem} \label{cbij}
 
 $$|S^{C}_k(n)|=|M^{C}_{n-k}(n)|, \mbox{ } k \in \{0\}\cup [n] .$$ 
 \end{theorem}
 
We prove Theorem \ref{cbij} by providing a bijection between the sets $S^{C}_k(n)$ and $M^{C}_{n-k}(n)$, $k \in \{0\}\cup [n]$. Before proceeding to the proof of  Theorem \ref{cbij}   we partition the sets  $S^{C}_k(n)$ and $M^{C}_{n-k}(n)$, $k \in \{0\}\cup [n]$, further. 
 
 Partition  $S^{C}_k(n)$, $k  \in \{0\}\cup [n]$, according to the $k$ distinct nonzero absolute values of the numbers appearing in the  sequence  and the number of times they appear. If $$\{|a_1|, |a_2|, \ldots, |a_n|\} \backslash \{0\}=\{c_1<c_2<\cdots<c_k\}$$ and $c_i$ appears $o_i$ times in $(|a_1|, |a_2|, \ldots, |a_n|)$, $ i \in [k]$, let $${S^{C}_k}^{\c, \o}(n)=\{(a_1, a_2, \ldots, a_n) \in S^{C}_k(n) |  
 \{\{|a_1|, |a_2|, \ldots, |a_n|\}\}=\cup_{i=1}^{k}\cup_{j=1}^{o_i}\{\{c_i\}\} \cup_{i=1}^{n-\sum_{j=1}^k o_j} \{\{0\}\}  \mbox{ } \},$$
 
 \noindent where $\c=(c_1< \ldots< c_k)$, $\o=(o_1, \ldots, o_k)$, $o_i>0$, for $i\in [k]$, and $\sum_{i=1}^k o_i\leq n$.
 
 Given an antichain $\a\in M^{C}_{n-k}(n)$, $k \in \{0\}\cup [n]$, it naturally corresponds to  a nonnesting $C_n$-partition $\pi_{\aa}$  of $[\pm n]$ with $k$ pairs of nonzero blocks. Let $(S_{\pi_{\a}}, g_{\pi_{\a}})$ be the pair of $k$-set and function corresponding to $\pi_{\a}$ under the bijection described in Theorem \ref{1}. Let  $$S_{\pi_{\a}}=\{c_1< \ldots< c_k\} \text{ and } o_i=g_{\pi_{\a}}(c_i), i \in [k].$$ Denote   $c(\aa)=(c_1< \ldots< c_k)$ and  $o(\aa)=(o_1, \ldots, o_k).$
 
  Partition the multiset $M^{C}_{n-k}(n)$, $k \in \{0\} \cup  [n]$, according to $\c=(c_1< \ldots< c_k)$, $\o=(o_1, \ldots, o_k)$, $o_i>0$, for $i\in [k]$, and $\sum_{i=1}^k o_i\leq n$, as described above. Let 

$${M^{C}_{n-k}}^{\c, \o}(n)=\{\{\aa \in M^{C}_{n-k}(n)|  
c(\aa)=\c, o(\aa)=\o \}\},$$
 
 \noindent where $\c=(c_1< \ldots< c_k)$, $\o=(o_1, \ldots, o_k)$, $o_i>0$, for $i\in [k]$, and $\sum_{i=1}^k o_i\leq n$.

\begin{lemma} \label{cun}
The vectors $c(\aa)=\c$ and $o(\aa)=\o$, where  $\c=(c_1< \ldots< c_k)$, $\o=(o_1, \ldots, o_k)$, $k \in \{0\}\cup  [n]$, $o_i>0$, for $i\in [k]$, $\sum_{i=1}^k o_i\leq n$,  uniquely determine the antichain $\aa$.
\end{lemma}
 
 \proof Lemma \ref{cun} follows readily since Theorem \ref{1} establishes a bijection.
 \qed

 \begin{theorem} \label{cco}
 $$|{S^{C}_k}^{\c, \o}(n)|=|{M^{C}_{n-k}}^{\c, \o}(n)|={n \choose {o_1, \ldots, o_k, n-\sum_{j=1}^k o_j}}2^{\sum_{j=1}^k o_j}, $$
 
 \noindent where $k \in \{0\} \cup [n]$, $\c=(c_1< \ldots< c_k)$, $\o=(o_1, \ldots, o_k)$, $o_i>0$, for $i\in [k]$, and $\sum_{i=1}^k o_i\leq n$.
 \end{theorem}
  
    \proof A bijective proof  can be given using Theorem \ref{1} and the ideas of Theorem \ref{3}. The enumeration is in Lemma \ref{2}.  Note that arrangements do not enter the proof.  \qed
    \medskip

 \noindent \textit{Proof of Theorem \ref{cbij}.} Straightforward corollary of Theorem \ref{cco}, since $$S^{C}_{k}(n)=\sum_{\c, \o} {S^{C}_{k}}^{\c, \o}(n)=\sum_{\c, \o} {M^{C}_{n-k}}^{\c, \o}(n)=M^{C}_{n-k}(n),$$  where $\c=(c_1< \ldots< c_k)$, $\o=(o_1, \ldots, o_k)$, $k \in [n]$, $o_i>0$, for $i\in [k]$, $\sum_{i=1}^k o_i\leq n$.
 \qed

 \begin{corollary}\label{cor:bb}
 
$$ \sum_{w \in \S_n^B} j(Q^{C}_w)=(2n+1)^{n}.$$
 \end{corollary}
 
 \proof Theorems \ref{cco} and \ref{cbij} extend to a bijection between $$\mathcal{M}^{C}(n)=\cup_{k=0}^n \cup_{\c, \o}{M^{C}_{n-k}}^{\c, \o}(n) \mbox{ and }\A^{C}(n)=\cup_{k=0}^n \cup_{\c, \o}{S^{C}_{k}}^{\c, \o}(n),$$ the cardinalities of which are  $ \sum_{w \in \S_n^B} j(Q^{C}_w)$ and $(2n+1)^{n}$, respectively.
 
\qed

  \section*{Acknowledgement}
  I would like to thank Richard Stanley for the beautiful problems he  poses in his classes, and which served as an inspiration for this paper. I would also like to thank Drew Armstrong for the many thoughtful suggestions and references he provided.

\end{document}